\documentclass[10pt]{article}

\usepackage{a4wide}
\usepackage{amssymb}
\usepackage{amsfonts}
\usepackage{amsmath}
\input xy
\xyoption{arrow} \xyoption{matrix}

\date{}

\newtheorem{proposition}{Proposition}[section]
\newtheorem{theorem}[proposition]{Theorem}
\newtheorem{lemma}[proposition]{Lemma}

\newtheorem{corollary}[proposition]{Corollary}

\def\der{\partial }

\def\nFM0{{\nu }_{F,M_0}}
\def\nFN0{{\nu }_{F,N_0}}
\def\nGN0{{\nu }_{G,N_0}}

\def\N0{ {\bf N}_0 }

\def\ra{\rightarrow}

\def\Xpm{X^{\pm }}

\def\s{\sigma}

\def\l1{{\lambda}_1}

\def\a{\alpha}
\def\a0{ {\alpha }_0}
\def\a1{ {\alpha }_1}

\def\l{\lambda}

\def\nFGM0{{\nu }_{F,G,M_0}}

\def\nFN0{{\nu}_{F,N_0}}

\def\sm{{\sigma}^m}

\def\sm1{{\sigma}^{-1}}

\def\smtp1{{\sigma}^{-t+1}}

\def\S1{S^{-1}}

\def\Xpm1{X^{\pm 1}_1}

\def\sPM1{{\sigma }^{\pm 1}}
\def\sMP1{{\sigma }^{\mp 1 }}

\def\d{\delta}

\def\di{{\rm d.ind}}

\def\L{\Lambda}

\def\CD{{\cal D}}

\def\Ytm1{Y^{t-1}}
\def\Yim1{Y^{i-1}}

\def\CG{{\cal G}}
\def\CH{{\cal H}}

\def\Aut{{\rm Aut}}

\def\Der{{\rm Der }}

\def\CJ{ {\cal J}}

\def\SL2Z{ {\rm SL}_2({\bf Z}) }

\def\Gp1{ G^{1 , 1 } }
\def\P11{ P^{-1 , 1 } }
\def\Pp1{ P^{1 , 1 } }

\def\nCLsr{{}^\nu\kern-2pt {\cal L}^{\sigma , \rho  }}
\def\nP{{}^\nu \kern-2pt P}
\def\nL{{}^\nu\kern-2pt L}
\def\nLL{{}^\nu\kern-2pt \Lambda}
\def\nPsr{{}^\nu\kern-2pt P^{\sigma , \rho  }}
\def\nLsr{{}^\nu\kern-2pt L^{\sigma , \rho  }}
\def\nuCL{{}^\nu\kern-2pt  {\cal L}}
\def\nCLsr{{}^\nu\kern-2pt {\cal L}^{\sigma , \rho  }}
\def\nCL1m{{}^\nu\kern-2pt {\cal L}^{-1 , 1  }}
\def\x1nu{x^\frac{1}{\nu}}
\def\xm1nu{x^{-\frac{1}{\nu}}}

\def\ra{\rightarrow }

\def\CB{{\cal B}}

\def\CH{ {\cal H}}

\def\nAM0{{\nu }_{{\cal A},M_0}}
\def\nAN0{{\nu }_{{\cal A},N_0}}

\def\End{ {\rm End }}
\def\Der{ {\rm Der }}
\def\CJ{ {\cal J }}

\def\det{ {\rm det }}

\def\gm{\mathfrak{m}}

\def\GL{{\rm GL}}
\def\SL{{\rm SL}}

\def\di!{\frac{\der^i}{i!}}
\def\dik!{\frac{\der^k_i}{k!}}

\def\N{\mathbb{N}}

\def\0{\overline{0}}
\def\1{\overline{1}}

\def\Ln1{\L_{n,\overline{1}}}

\def\a1{a_{\overline{1}}}

\def\S{\Sigma}

\def\grad{{\rm grad}}

\def\vn1{\overrightarrow{n-1}}

\def\Sh{{\rm Sh}}

\def\mS{\mathbb{S}}
\def\mJ{\mathbb{J}}
\def\mI{\mathbb{I}}

\def\K1{{\rm K}_1}

\def\hmI1{\widehat{\mI_1}}
\def\tmI1{\widetilde{\mI_1}}
\def\tmJ1{\widetilde{\mJ_1}}
\def\hB1{\widehat{B_1}}
\def\hCB1{\widehat{\CB_1}}

\def\Fix{{\rm Fix}}

\def\mJ{\mathbb{J}}

\def\AutKalg{ {\rm Aut_{K-{\rm alg}}}}

\def\mmW1{\mathbb{W}_1}

\def\divn0{\mathfrak{div}_n^0}
\def\div0mu{\mathfrak{div}_{n, [\mu ]}^0}
\def\din0{\mathfrak{di}_n^0}
\def\ivn0{\mathfrak{iv}_n^0}
\def\divnc{\mathfrak{div}_n^c}
\def\divv{{\rm div}}
\def\mfGn{\mathbf{G}_n}
\def\mfGnc{\mathbf{G}_n^c}

\def\AutKalgc{ {\rm Aut_{K-{\rm alg},c}}}

\def\gsn{\mathfrak{s}_n}

\def\Divn0{\mathfrak{Div}_n^0}
\def\Din0{\mathfrak{Di}_n^0}
\def\Ivn0{\mathfrak{Iv}_n^0}
\def\Divnc{\mathfrak{Div}_n^c}
\def\mfGhn{\widehat{\mathbf{G}}_n}
\def\mfGhnc{\widehat{\mathbf{G}}_n^c}

\begin{document}

\author{V. V. \  Bavula  
}

\title{The groups of automorphisms of the Lie algebras of formally analytic vector fields with constant divergence}

\maketitle

\begin{abstract}
Let $S_n=K[[x_1, \ldots , x_n]]$ be the  algebra of power series  over a field $K$ of characteristic zero, $\mS_n^c$ be the group of continuous automorphisms of $S_n$ with constant Jacobian, and    $\Divnc$ be the Lie algebra of derivations of $S_n$ with constant divergence.
We prove that $\Aut_{{\rm Lie}}(\Divnc )=\Aut_{{\rm Lie}, c}(\Divnc )\simeq \mS_n^c$.

$\noindent $



\end{abstract}


\section{Introduction}

In this paper,  $K$ is a
field of characteristic zero and  $K^*$ is its group of units, and the following notation is fixed:
\begin{itemize}
\item $P_n:= K[x_1, \ldots , x_n] $ 
  is a polynomial algebra  and
$S_n:= K[[x_1, \ldots , x_n]]$ is the algebra of power series over $K$, $\gm := (x_1, \ldots , x_n)$,  and $S_n^*$ is the group of units of $S_n$,
 \item 
      $\mS_n:=\AutKalgc (S_n)$ is the group of continuous (with respect to the $\gm$-adic topology) automorphisms of  $S_n$   and $\mS_n^c:=\{ \s \in \mS_n\, | \, \CJ (\s ) \in K\}$ where $\CJ (\s )$ is the Jacobian of $\s$,

 \item $\der_1:=\frac{\der}{\der x_1}, \ldots , \der_n:=\frac{\der}{\der
x_n}$ are the partial derivatives ($K$-linear derivations) of
$S_n$,
\item    
$\gsn:=\Der_K(S_n) =\bigoplus_{i=1}^nS_n\der_i$ 
  is the Lie
algebras of $K$-derivations of $S_n$ where $[\der , \d ]:= \der \d -\d \der $, and $D_n:= \Der_K(P_n) =\bigoplus_{i=1}^n P_n\der_i$,
 \item $\CD_n:=\bigoplus_{i=1}^n K\der_i$,
 \item 
 $H_1:=x_1\der_1, \ldots , H_n:=x_n\der_n$,

 \item for a derivation $\der = \sum_{i=1}^n a_i\der_i\in \gsn$, $\divv (\der ) := \sum_{i=1}^n\frac{\der  a_i}{\der x_i}$ is the {\em divergence}  of $\der$,
\item $\divn0 :=\{ \der \in D_n \, | \, \divv (\der ) =0\}$  and $\Divn0 :=\{ \der \in \gsn \, | \, \divv (\der ) =0\}$ are  the Lie algebras  of polynomial, respectively, formally analytic  vector fields (derivations) with zero divergence,
 \item $\mfGn :=\Aut_{{\rm Lie}}(\divn0 )$ and $\mfGhn :=\Aut_{{\rm Lie}}(\Divn0 )$,
 \item $\divnc :=\{ \der \in D_n \, | \, \divv (\der ) \in K\}$ and $\Divnc :=\{ \der \in \gsn \, | \, \divv (\der ) \in K\}$ are the Lie algebras  of polynomial, respectively, formally analytic vector fields (derivations) with constant  divergence,

 \item $\mfGnc :=\Aut_{{\rm Lie}}(\divnc )$ and $\mfGhnc :=\Aut_{{\rm Lie}}(\Divnc )$,
    \end{itemize}

{\bf The groups of automorphisms of the Lie algebras $\divn0$ and $\divnc$}. Let $\Sh_1:=\{ s_\mu \in \AutKalg (K[x])\, | \, s_\mu (x) = x+\mu, \; \mu \in K\}$.
\begin{theorem}\label{16Mar13}
\cite{Rudakov-1986}, \cite{Bav-Aut-Div} $\mfGn \simeq \begin{cases}
G_1/\Sh_1\simeq K^*& \text{if }n=1,\\
G_n& \text{if }n\geq 2.\\
\end{cases}$
\end{theorem}
Theorem \ref{16Mar13} was announced  in \cite{Rudakov-1986}  where a  sketch of the proof is given based on a study of certain Lie subalgebras of $\divn0$ of finite codimension. The proof in \cite{Bav-Aut-Div}  is based on completely different ideas.
\begin{theorem}\label{A16Mar13}
\cite{Bav-Aut-Div} $\mfGnc \simeq  G_n$.
\end{theorem}

{\bf The groups of automorphisms of the Lie algebras $\Divn0$ and $\Divnc$}.

\begin{theorem}\label{FD16Mar13}
\cite {Rudakov-1969}, \cite{Rudakov-1986} $\mfGhn \simeq \mS_n^c $ for $n\geq 2$.
\end{theorem}


The aim of the paper is to prove the following theorem.

\begin{theorem}\label{FDA16Mar13}
 $\mfGhnc \simeq \begin{cases}
G_1& \text{if }n=1,\\
 \mS_n^c& \text{if }n\geq 2.\\
\end{cases}$
\end{theorem}

{\em Proof}. For $n=1$, $\mathfrak{Div}_1^c = K\der_1\oplus KH_1= \mathfrak{div}_1^c$ and so
$\widehat{\mathbf{G}}_1^c=\mathbf{G}_1^c=G_1$, by Theorem \ref{A16Mar13}. So, let $n\geq 2$.

$\noindent $

(i) $\mS_n^c \subseteq \mfGhnc$ via the group monomorphism (Theorem \ref{FD16Mar13} and Theorem \ref{FB16Mar13}):
$$\mS_n^c\ra \mfGhnc, \;\;  \s \mapsto \s : \der \mapsto \s (\der ):=\s \der \s^{-1}.$$


(ii) $\Divn0 = [\Divnc , \Divnc ]$: The equality follows from the fact that $\Divn0$ is a simple Lie algebra which is an ideal of the Lie algebra $\Divnc$ and $\Divnc = \Divn0 \oplus KH_1$.

$\noindent $

(iii) The short exact sequence of group homomorphisms
$$ 1\ra F:=\Fix_{\mfGhnc}(\Divn0 ) \ra \mfGhnc
\stackrel{{\rm res}}{\ra}\mfGhn\ra 1$$
is exact  (by (i) and Theorem \ref{FD16Mar13}) where res $:\s \mapsto \s|_{\Divn0}$ is the restriction map, see (ii).

 $\noindent $

 (iv)  Since $\mfGhn = \mS_n^c$ (Theorem \ref{FD16Mar13}) and $\mS_n^c\subseteq \mfGhnc$ (by (i)), the short exact sequence splits
\begin{equation}\label{FDG=GF}
 \mfGhnc \simeq \mfGhn\ltimes F.
\end{equation}


(v) $ F=\{ e\}$ (Lemma \ref{Fa18Mar13}). Therefore, $\mfGhnc \simeq  \mS_n^c$.   $\Box $

$\noindent $

The Lie algebra $\Divnc$ is a topological Lie algebra with respect to the $\gm$-adic topology, i.e., the set $\{ \gm^i \Divnc\}_{i\in \N}$ is a base of open neighbourhoods of zero. Let
 $\widehat{\mathbf{G}}_{n,top}^c$ be group of automorphisms of the topological Lie algebra $\mathfrak{Div}_n^c$. Clearly, $\widehat{\mathbf{G}}_{n,top}^c\subseteq\mfGhnc$. The inverse inclusion follows from Theorem \ref{FDA16Mar13}.

\begin{corollary}\label{aFDA16Mar13}
$\widehat{\mathbf{G}}_{n,top}^c=\mfGhnc$.
\end{corollary}

$\noindent $

{\bf The group $\mS_n$}.
Every continuous automorphism $\s \in \mS_n$  is uniquely determined by the elements
$$x_1':=\s (x_1), \; \ldots , \; x_n':=\s (x_n)$$
that necessarily  (by the continuity of $\s $) belong to the maximal ideal $\gm$ of the algebra $S_n$, and for all  series $f=f(x_1, \ldots , x_n)\in S_n$, $\s (f) = f(x_1', \ldots , x_n')$.
Let $M_n(S_n)$ be the algebra of $n\times n$ matrices over  $S_n$. The matrix  $J(\s) := (J(\s )_{ij}) \in M_n(S_n)$, where $J(\s )_{ij} =\frac{\der x_j'}{\der x_i}$,   is called the {\em Jacobian matrix} of  $\s$ and its determinant $\CJ (\s ) :=\det \, J(\s)$ is called the {\em Jacobian} of $\s$. So, the $j$'th column of $J(\s )$ is the {\em gradient} $\grad \, x_j':=(\frac{\der x_j'}{\der x_1}, \ldots , \frac{\der x_j'}{\der x_n})^T$  of the series $x_j'$. Then the derivations
$$\der_1':= \s \der_1\s^{-1}, \; \ldots , \; \der_n':= \s\der_n\s^{-1}$$ are the partial derivatives of $S_n$ with respect to the variables $x_1', \ldots , x_n'$,
\begin{equation}\label{ddp=dxi}
\der_1'=\frac{\der}{\der x_1'}, \; \ldots , \; \der_n'=\frac{\der}{\der x_n'}.
\end{equation}
Every derivation $\der \in \gsn$ is a unique sum $\der = \sum_{i=1}^n a_i\der_i$ where $a_i = \der *x_i\in S_n$. Let  $\der := (\der_1, \ldots , \der_n)^T$ and $ \der' := (\der_1', \ldots , \der_n')^T$ where $T$ stands for the transposition. Then
\begin{equation}\label{dp=Jnd}
\der'=J(\s )^{-1}\der , \;\; {\rm i.e.}\;\; \der_i'=\sum_{j=1}^n (J(\s )^{-1})_{ij} \der_j\;\; {\rm for }\;\; i=1, \ldots , n.
\end{equation}
In more detail, if $\der'=A\der $ where $A= (a_{ij})\in M_n(S_n)$, i.e. $\der_i=\sum_{j=1}^n a_{ij}\der_j$. Then for all $i,j=1, \ldots , n$,
$$\d_{ij}= \der_i'*x_j'=\sum_{k=1}^na_{ik}\frac{\der x_j'}{\der x_k}$$
where $\d_{ij}$ is the Kronecker delta function. The equalities above can be written in the matrix form as  $AJ(\s) = 1$ where $1$ is the identity matrix. Therefore, $A= J(\s )^{-1}$.

For all $\s, \tau \in \mS_n$,
\begin{equation}\label{SJst=JsJ}
J(\s \tau ) = J(\s )\cdot  \s (J(\tau )).
\end{equation}
By taking  the determinants of both sides of (\ref{SJst=JsJ}), we have a similar equality of the Jacobians: for all $\s , \tau \in \mS_n$.
\begin{equation}\label{SJst=JsJ1}
\CJ (\s \tau ) = \CJ(\s )\cdot  \s (\CJ(\tau )).
\end{equation}
By putting $\tau = \s^{-1}$ in (\ref{SJst=JsJ}) and (\ref{SJst=JsJ1}) we see that $\ J(\s ) \in \GL_n(S_n)$,   $\CJ (\s ) \in S_n^*$, and
\begin{equation}\label{SJsm}
J(\s^{-1}) = \s^{-1} (J(\s )^{-1}),
\end{equation}
\begin{equation}\label{SJsm1}
\CJ (\s^{-1} ) =  \s^{-1} (\CJ(\s )^{-1}).
\end{equation}
\begin{eqnarray*}
 \mS_n &=& \{ \s \in \End_{K-{\rm alg},c}(S_n)\, | \, \CJ (\s ) \in S_n^* \} \\
 &=& \{ \s \in \End_{K-{\rm alg},c}(S_n)\, | \, \s (x) = Ax+\cdots , A= (a_{ij}) \in \GL_n(K)\} ,
\end{eqnarray*}
that is $\s (x_i) = \sum_{j=1}^n a_{ij} x_j+\cdots $ where the three dots mean smaller terms ($\cdots \in  \gm^2$).

$\noindent $

\begin{lemma}\label{a15Oct13}
For all $\s \in \mS_n^c$, $$ \sum_{j=1}^n \der_j*(J(\s )^{-1})_{ij}=0 \;\; {\rm for } \;\; i=1, \ldots , n.$$
\end{lemma}

{\it Proof}. By (\ref{dp=Jnd}), $\der_i' = \sum_{i=1}^n (J(\s )^{-1})_{ij}\der_j$. By Theorem \ref{FD16Mar13}, we have the result.
$\Box $

$\noindent $

{\bf The divergence commutes with automorphisms $\mS_n^c$}. The following theorem shows that the divergence commute with automorphisms $\mS_n^c$, i.e. the divergence map $\divv : \gsn \ra S_n$ is an  $\mS^c_n$-module homomorphism.

\begin{theorem}\label{FB16Mar13}
For all $\s \in \mS^c_n$ and $\der \in \gsn$,
$$ \divv (\s (\der ))= \s (\divv (\der )).$$
\end{theorem}

{\it Proof}. Let $\der = \sum_{i=1}^n a_i\der_i$ where $a_i\in S_n$. Then
$ \der' = \s \der \s^{-1} =\sum_{i=1}^n \s (a_i)  \der_i'$ where, by (\ref{dp=Jnd}), $\der_i' = \sum_{j=1}^n (J(\s )^{-1})_{ij}\der_j$. Now, by Lemma \ref{a15Oct13},
\begin{eqnarray*}
 \divv (\der') &=& \sum_{i,j=1}^n \der_j*((J(\s)^{-1})_{ij}\s (a_i)) = \sum_{i=1}^n (\sum_{j=1}^n \der_j *(J(\s )^{-1})_{ij})\cdot \s (a_i))+ \sum_{i=1}^n \sum_{j=1}^n (J(\s)^{-1})_{ij} \der_j*\s (a_i) \\
 &=&  \sum_{i=1}^n \der_i'*\s (a_i) = \sum_{i=1}^n \s \der_i \s^{-1}\s (a_i) = \s ( \sum_{i=1}^n \der_i (a_i)) = \s (\divv (\der )). \;\; \Box
\end{eqnarray*}

\begin{lemma}\label{Fa18Mar13}
 $ \Fix_{\mfGhnc}(\Divn0 ) =\{ e\}$ for $n\geq 2$.
\end{lemma}

{\it Proof}.  Let $\s \in F:= \Fix_{\mfGhnc}(\Divn0 )$, $H_1':= \s (H_1) , \ldots , H_n':= \s (H_n)$. Since $\Divnc= \Divn0\oplus KH_i$, $i=1, \ldots , n$, it suffices to show that $ \s (H_i) = H_i$ for $i=1, \ldots , n$. For $i\neq j$, $\s (H_i-H_j) = H_i-H_j$, and so $d:= H_i'-H_i= H_j'-H_j$. For all $i=1, \ldots , n$,
$$ [\der_i , d] = \s ([\der_i, H_i]) - [\der_i, H_i]= \s (\der_i ) - \der_i = \der_i - \der_i =0.$$ So, $d\in C_{\Divnc}(\CD_n ) = \CD_n$ (since $C_{\gsn}(\CD_n ) = \CD_n$) and $d= \sum_{i=1}^n \l_i\der_i$ for some $\l_i\in K$ where $C_\CG (\CH ):= \{ g\in \CG \, | \, [g, \CH ]=0\}$ is the centralizer of a subset $\CH$ of a Lie algebra $\CG$. The elements $H_1'=H_1+d, \ldots , H_n'=H_n+d$ commute hence $d=0$. Therefore, $\s = e$.  $\Box $

$\noindent $

$${\bf Acknowledgements}$$

 The work is partly supported by  the Royal Society  and EPSRC.

\small{

Department of Pure Mathematics

University of Sheffield

Hicks Building

Sheffield S3 7RH

UK

email: v.bavula@sheffield.ac.uk}

\end{document}